\documentclass[a4paper,10pt]{article}%
\usepackage[english]{babel}
\usepackage[utf8]{inputenc}

\usepackage{amsmath}
\usepackage{amsfonts}
\usepackage{amssymb}
\usepackage{amsthm}
\usepackage{graphicx}
\usepackage{pdfpages}

\newcommand{\mytext}[1]{ \: \textrm{#1} \: }
\newcommand{\mysetdescr}[2]{\left\{ #1 \: \left| \: #2 \right. \right\} }
\newcommand{\mydownarrow}{{\downarrow \,}}

\newcommand{\myN}{\mathbb{N}}
\newcommand{\myNk}[1]{\underline {#1}}

\newcommand{\Nk}{\myNk{k}}
\newcommand{\Nn}{\myNk{n}}

\newcommand{\Nke}{\myNk{k-1}}

\def\D{{\cal D}}
\def\E{{\cal E}}
\def\H{{\cal H}}
\def\J{{\cal J}}
\def\P{{\cal P}}

\newcommand{\mf}[1]{\mathfrak{ #1 }}
\newcommand{\fp}{\mf{p}}
\newcommand{\fP}{\mf{P}}

\def\BP{\begin{proof}}
\def\EP{\end{proof}}

\newcommand{\equivA}{ \; \equiv \; }
\newcommand{\simeqA}{ \; \simeq \; }
\newcommand{\eqA}{ \; = \; }
\newcommand{\inA}{ \in \; }
\newcommand{\subseteqA}{ \subseteq \; }
\newcommand{\capA}{ \cap \; }
\newcommand{\cupA}{ \cup \; }

\newcommand{\myHCC}[2]{ \H(C_{#1}, C_{#2} ) }
\newcommand{\myHCCC}[3]{ \H( \myHCC{#1}{#2}, C_{#3} ) }
\newcommand{\Hzk}{\myHCC{2}{k}}

\newcommand{\Dzk}{\D( \Hzk ) }

\newcommand{\smda}[1]{ \; \setminus \! \downarrow_{#1} \! }
\newcommand{\smdaR}{ \smda{R} }
\newcommand{\smdaPp}{ \smda{P^+} }

\newcommand{\CC}[2]{ C_{#1} \times C_{#2} }
\newcommand{\CCnk}{ \CC{n}{k} }

\newcommand{\WC}[1]{ W \times C_{#1} }
\newcommand{\WCk}{ \WC{k} }

\newcommand{\WCkm}{ \WC{k-1} }

\newcommand{\LC}[1]{ \Lambda \times C_{#1} }
\newcommand{\LCk}{ \LC{k} }

\newcommand{\LCkm}{ \LC{k-1} }

\newcommand{\CCdrei}[1]{ C_3 \times C_3 \times C_{#1} }
\newcommand{\CCdreik}{ \CCdrei{k} }

\newcommand{\LOC}[1]{ \lozenge \times C_{#1} }
\newcommand{\LOCk}{ \LOC{k} }

\newcommand{\LOCkm}{ \LOC{k-1} }

\newcommand{\Mell}{ {\{ \ell \}} }
\newcommand{\Merr}{ {\{ r \}} }
\newcommand{\Mellerr}{ {\{ \ell, r \}} }

\newcommand{\Mellb}{ {\{ \bot, \ell \}} }
\newcommand{\Merrb}{ {\{ \bot, r \}} }
\newcommand{\Mellerrb}{ {\{ \bot, \ell, r \}} }

\begin{document}

\theoremstyle{plain}
\newtheorem{condition}{Condition}
\newtheorem{theorem}{Theorem}
\newtheorem{definition}{Definition}
\newtheorem{corollary}{Corollary}
\newtheorem{lemma}{Lemma}
\newtheorem{proposition}{Proposition}

\title{\bf A recursive approach for the enumeration of the homomorphisms from a poset $P$ to the chain $C_3$}
\author{\sc Frank a Campo}
\date{\small Viersen, Germany\\
{\sf acampo.frank@gmail.com} \\
April 2021}

\maketitle

\begin{abstract}
\noindent Let $\H(P,C_3)$ be the set of order homomorphisms from a poset $P$ to the chain $C_3 = 1 < 2 < 3$. We develop a recursive approach for the calculation of the cardinality of $\H(P,C_3)$, and we apply it on several types of posets, including $P = C_3 \times C_3 \times C_k$ and $P = \H(C_k, C_3)$; for the latter poset $P$, we derive a direct formula for $\# \H ( P, C_3 )$.
\newline

\noindent{\bf Mathematics Subject Classification:}\\
Primary: 06A07. Secondary: 06A06.\\[2mm]
{\bf Key words:} poset, homomorphism, chain.
\end{abstract}

\section{Introduction} \label{sec_introduction}

Let $C_n$ be the set $\{ 1, \ldots , n \}$ equipped with the natural order. The number of homomorphisms from a poset $P$ to $C_n$ is the value at $x = n$ of the {\em order polynomial} $\Omega_P(x)$ of $P$, introduced by Stanley in the early seventies \cite{Stanley_1970,Stanley_1971_Diss,Stanley_1972}. Due to \cite[Theorem 2]{Stanley_1970}, the order polynomial can be written as
\begin{align} \label{formel_omegaPx}
\Omega_P(x) & \eqA \sum_{d=0}^{\# P -1} w_P(d) \binom{\# P + x -1-d}{\# P },
\end{align}
where $w_P(d)$ is the number of linear extensions of $P$ (regarded as permutations of a natural labeling of $P$) with exactly $d$ descents.

The research about the order polynomial focuses on relating it to other structures and polynomials. The following incomplete list of topics and references highlights the variety of subjects. The connections between the order polynomial, the chromatic polynomial in graph theory, and the Erhart polynomial of an order polytope have already been seen by Stanley \cite{Stanley_1970,Stanley_1984}. Edelman and Klingsberg \cite{Edelman_Klingsberg_1982} use the lattice of sub-posets of a poset in order to prove identities for the order polynomial, and Wagner \cite{Wagner_1992} asks for the zeros of another polynomial related to it. Hamaker et al.\ \cite{Hamaker_etal_2020} and Browning et al.\ \cite{Browning_etal_2018} study posets with identical order polynomial, and Jochemko \cite{Jochemko_2014} connects the concepts of the order polynomial and P\'{o}lya's enumeration theorem. The general formula presented by Thomas \cite{Thomas_2003} expresses $\Omega_P(n)$ by means of coefficients related to the order polytope of $P$; just as \eqref{formel_omegaPx}, Thomas' formula provides structural insight, but is not intended to be a tool for practical calculation.

Direct attempts to really calculate the order polynomial or its values are restricted in most cases to examples and exercises \cite[Abschnitt III.4]{Aigner_1975}, \cite[Section 3.15, in part.\ Ex.\ 3.62, 3.66]{Stanley_2012}. For posets with an appropriate structure, the order polynomial has a product formula and is thus (at least in principle) easier to calculate and to evaluate; a survey with references is presented by Hopkins \cite{Hopkins_2020}. The reason for the reserve to calculate $\Omega_P(x)$ or its values is the complexity of the task: Brightwell and Winkler \cite{Brightwell_Winkler_1991} showed that the calculation of $\Omega_P( \# P )$ (the number of linear extensions of $P$) is $\#P$-complete.

This number of linear extensions is according to \cite[p.\ 258]{Stanley_2012} probably the single most useful number for measuring the ``complexity'' of a poset. An own branch of mathematics has developed around it, and presumably every reader has already been in touch with it. However, there is also some interest to know the values $\Omega_P(n)$ for {\em small} integers $n$. $\Omega_P(1) = 1$ is trivial, and with $\D(P)$ being the down-set lattice of $P$, the equation $\Omega_P(2) = \# \D(P)$ belongs to the basics of order theory, and the enumeration of $\D(P)$ is a standard task which often can be done manually. The numbers $h( P ) \equiv \Omega_P(3) = \# \H(P,C_3)$ are considerably more complicated to determine, and we present in this article a recursive method for their effective calculation.

After recalling common terms and notation in Section \ref{sec_notation}, we develop our approach in Section \ref{sec_recursion}. Our first starting point is the well-known isomorphism $\H( P, C_3 ) \simeq \H( C_2, \D(P) )$ yielding
\begin{align*}
h(P) & \eqA \sum_{D \in \D(P)} \# \mydownarrow_{\D(P)} D,
\end{align*}
a summation formula widely used, e.g., in the computation of Dedekind numbers \cite{aCampo_2018}. The second  starting point is the flexible concept of the {\em generalized vertical sum} of posets introduced by the author and Ern\'{e} \cite{aCampo_2019,aCampo_Erne_2019}. For a generalized vertical sum $R$ of two posets $P$ and $Q$, it has been shown \cite{aCampo_Erne_2019} that the down-set lattice $\D(R)$ of $R$ is the disjoint union of certain sets $\J_T(R)$, with $T$ running through $\D(Q)$; we recall this result in Theorem \ref{theo_genvertsum}. With
\begin{align*}
a_T(R) & \equivA \sum_{D \in \J_T(R)} \# \mydownarrow_{\D(R)} D \quad \mytext{for all } T \in \D(Q),
\end{align*}
we thus get
\begin{align*}
h(R) & \eqA  \sum_{T \in \D(Q)} a_T(R).
\end{align*}
We now postulate without loss of generality that the posets $P$ and $Q$ are linked in $R$ by two sub-posets $S^- \subseteq P$ and $S^+ \subseteq Q$ and a mapping $\sigma$ from $\D(S^+)$ to subsets of $S^-$ which is compatible with the structure of $R$. The sub-posets $S^+$ and $S^-$ give raise to several interrelated generalized vertical sums, and their respective sets $\J_{T'}(R')$ are linked in Lemma \ref{lemma_beta_tau} by isomorphisms. Based on these results, a recursive formula for the coefficients $a_T(R)$ with $T \in \D(S^+)$ is derived in Theorem \ref{theo_formel_a}. It refers to proper sub-posets of $R$ only and offers thus the possibility to calculate $h(R)$ recursively. The approach is very flexible and can be set up in different ways for a given poset $R$. It even gives raise to new ways to calculate Dedekind  numbers, because $h( C_2^k )$ is the $(k+1)$\textsuperscript{th} Dedekind number.

For posets $R$ fitting well to the structure of the recursion, $h(R)$ can be calculated manually. We do so in Section \ref{sec_application}. In Section \ref{subsec_LambdaDiamond}, we work with posets $R = W \times C_k$, and we calculate $h( W \times C_k)$ for several posets $W$, including the chain $C_n$, the poset $\Lambda$ with $\Lambda$-shaped diagram, the diamond $C_2 \times C_2$, and the Noughts and Crosses grid $C_3 \times C_3$. In Section \ref{subsec_HC2Ck}, we treat the posets $R = \H(C_2,C_k) \simeq \H( C_{k-1}, C_3 )$ and derive a closed formula for $h(R)$.

Besides of its mathematical interest, the number $h(R)$ has some relevance in other areas of science, too, e.g., for {\em ensemble based systems} in machine learning  \cite{Polikar_2006}, also known as {\em multiple classifier systems, committee of classifiers}, or {\em mixture of
experts}. Assume that $k$ experts (humans, robots, recognition systems, ...) have the task to rank objects on a scale with $v$ rank levels, e.g., ``stop'' and ``go'' for $v=2$, or ``negative'', ``neutral'', ``positive'' for $v=3$. The judgement results in a point $r \in C_v^k$, and now a summary value $s(r) \in C_v$ has to be assigned to $r$ by a {\em summary rule} $s$. Because a better ranking of the experts cannot downgrade the summary value, the possible summary rules are the elements of $\H( C_v^k, C_v)$, and fundamental questions ask for their number, classification etc. For $v = 2$, the summary rules are the monotone Boolean functions \cite{Kovalerchuk_1995,Liggins_Nebrich_2000}, and the figures $\# \H( C_2^k, C_2)$ are the Dedekind  numbers \cite{aCampo_2018} known up to $k = 8$. For $v = 3$, we deal with monotone {\em ternary} functions. The number $h( C_3^2 ) = 175$ can still be determined with paper and pencil, but already the calculation of $h( C_3^3 ) = 211250$ done in Section \ref{subsec_LambdaDiamond} is out of reach of manual calculation. Also Section \ref{subsec_HC2Ck} has a connection to summary rules, because $h( \H(C_2,C_k ) )$ is the number of {\em symmetric} summary rules.

\section{Basics and Notation} \label{sec_notation}

We are working with {\em finite posets}, thus ordered pairs $P = (X,\leq_P)$ consisting of a finite set $X$ (the {\em carrier} of $P$) and a {\em partial order relation} $\leq_P$ on $X$, i.e., a reflexive, antisymmetric, and transitive subset of $X  \times X$. Due to reflexivity, the {\em diagonal} $\Delta_X \equiv \mysetdescr{(x,x)}{x\in X}$ is always a subset of $\leq_P$. As usual, we write $x \leq_P y$ for $(x,y) \inA \leq_P$.

We say that $y \in P$ is {\em covered} by $x \in P$, iff $y \not= x$ and $y \leq_P x$ without any additional point between them: $y \leq_P z \leq_P x \Rightarrow z \in \{ x, y \}$ for all $z \in P$. 

Two elementary posets can be defined on any set $X$: The {\em antichain} $(X, \Delta_X)$ and the {\em chain} which is up to isomorphism characterized by $x \leq_P y$ or $y \leq_P x$ for all $x, y \in X$. For a finite set $X$ with $k \equiv \# X$, we write $A_k$ for the antichain on $X$ and $C_k$ for the chain on $X$. In what follows, $C_k$ is always the set $\{ 1 , \ldots , k \}$ equipped with the natural order.

A poset $Q = (Y, \leq_Q)$ is called a {\em sub-poset} of $P$ iff $Y \subseteq X$ and $\leq_Q \; \subseteqA \leq_P$, and for $Y \subseteq X$, the {\em induced poset} $P \vert_Y$ is defined as $\left( Y, \leq_P \capA (Y \times Y) \right)$; however, we write $P \setminus Y$ instead of $P \vert_{X \setminus Y}$.

Given two posets $P = (X, \leq_P)$ and $Q = (Y, \leq_Q )$, we can construct new posets with them. $P \times Q$ is the poset with carrier $X \times Y$ and the component-wise defined partial order relation. If $X$ and $Y$ are disjoint, the {\em direct sum} $P + Q$ and the {\em ordinal sum} $P \oplus Q$ are posets on $X \cup Y$ with partial order relations
\begin{align*}
\leq_{P+Q} & \equivA \leq_P \cup \leq_Q, \\
\leq_{P \oplus Q} & \equivA \leq_P \cup \leq_Q \cupA (X \times Y).
\end{align*}
The {\em generalized vertical sums} have been introduced by the author and Ern\'{e} \cite{aCampo_2019,aCampo_Erne_2019} as structures "in-between" direct sums and ordinal sums:
\begin{definition}[\cite{aCampo_2019,aCampo_Erne_2019}]  \label{def_genvertsum}
Let $P = (X, \leq_P), Q = (Y, \leq_Q)$ be posets with disjoint carriers $X$ and $Y$. A poset $R = ( X \cup Y, \leq_R )$ on $X \cup Y$ is called a {\em generalized vertical sum} of $P$ and $Q$ iff 
\begin{displaymath}
\leq_P \cup \leq_Q \quad \subseteq \quad \leq_R \quad \subseteq \;  \quad \leq_P \cup \leq_Q \cup \; (X \times Y).
\end{displaymath}
We call $P$ the {\em lower part} and $Q$ the {\em upper part} of $R$; ``generalized vertical sum'' is abbreviated as ``g.v.s.'' in what follows.
\end{definition}

Down-sets are one of the fundamental concepts in order theory. Given a poset $P = (X, \leq_P)$, a subset $D \subseteq X$ is called a {\em down-set} or {\em order ideal} of $P$ iff $x \in D$ holds for every $x \in X$ for which a $y \in D$ exists with $x \leq_P y$. For $B \subseteq X$ and $x \in X$, we define the down-sets created by $B$ and $x$ in $P$ by
\begin{align*}
\mydownarrow_P \; B & \equivA \mysetdescr{ a \in X }{ a \leq_P b \; \; \mytext{for a } b \in B }, \\
\mydownarrow_P \; x & \equivA \mydownarrow_P \; \{ x \}.
\end{align*}

The set of down-sets of $P$ is denoted by $\D(P)$. Together with set inclusion, $\D(P)$ is a partial order (even a lattice). For a down-set $D \in \D(P)$, the symbol $\mydownarrow_{\D(P)} D$ thus indicates the down-set created by $D$ in $\D(D)$:
\begin{align*}
\E_D(P) & \equivA \mydownarrow_{\D(P)} D \eqA \mysetdescr{ E \in \D(P) }{ E \subseteq D }.
\end{align*}

{\em Up-sets} are the duals of downsets: a subset $U \subseteq X$ is called an {\em up-set} or {\em order filter} of $P$ iff $x \in U$ holds for every $x \in X$ for which a $y \in U$ exists with $y \leq_P x$.

In order to make the line of thought more conclusive, we frequently identify a down-set of a poset with the poset induced by it, e.g., by calling $P$ a down-set of $P$.

A mapping $\xi : X \rightarrow Y$ is called an {\em (order) homomorphism} from a poset $P = (X, \leq_P)$ to a poset $Q = (Y, \leq_Q)$ iff $x \leq_P y$ implies $\xi(x) \leq_Q \xi(y)$ for all $x, y \in X$. The set of order homomorphisms from $P$ to $Q$ is denoted by $\H(P,Q)$. We make $\H(P,Q)$ being a poset by equipping it with the usual point-wise partial order $\leq_{\H(P,Q)}$ defined by
\begin{align*}
\xi \leq_{\H(P,Q)} \zeta & \equivA \xi(x) \leq_Q \zeta(x) \quad \mytext{for all } x \in X.
\end{align*}
``$\simeq$'' indicates isomorphism of posets.

From set theory, we use the following symbols:
\begin{align*}
\myNk{0} & \equiv  \emptyset, \\
\myNk{n} & \equiv  \{ 1, \ldots, n \} \mytext{for every} n \in \myN, \\
\myNk{n}_0 & \equiv \myNk{n} \cup \{ 0 \} \mytext{for every} n \in \myN_0,
\end{align*}
and for every set $X$, the symbol $\P(X)$ denotes the {\em power set} of $X$.



\section{The recursion} \label{sec_recursion}

For the determination of the cardinality of $\H(R,C_3)$, we start with the general isomorphism \cite[p.\ 4]{aCampo_2018}
\begin{align} \label{isom_HPDQ_HQDP}
\H(R, \D(Q)) & \; \simeq\; \H(Q, \D(R)),
\end{align}
yielding
\begin{align*}
\H( R, C_3 ) & \; \simeq \; \H( C_2, \D( R ) ).
\end{align*}
We can thus determine $\# \H( R, C_3 )$ by means of the summation formula
\begin{align} \label{hauptgleichung}
h(R) & \eqA \sum_{D \in \D(R)} \# \mydownarrow_{\D(R)} D \eqA \sum_{D \in \D(R)} \# \E_D(R) .
\end{align}
This formula is well-known and has widely been used in the computation of Dedekind numbers \cite{aCampo_2018}.

Fundamental for our apporoach is the following theorem describing the down-set lattice of a generalized vertical sum:
\begin{theorem}[\cite{aCampo_Erne_2019}, Theorem 3.3]  \label{theo_genvertsum}
Let $P = (X, \leq_P), Q = (Y, \leq_Q)$ be posets with disjoint carriers $X$ and $Y$ and let $R = ( X \cup Y, \leq_R )$ be a g.v.s.\ with lower part $P$ and upper part $Q$. Then the down-set lattice of $R$ is given by the following disjoint union:
\begin{align} \label{formel_downsets}
\D(R) & = \; \bigcup_{T \in \D(Q)}
\mysetdescr{ D \; \cup \downarrow_R T }{ D \in \D( P \smdaR T) }.
\end{align}
\end{theorem}
In what follows, the symbols $P$, $Q$, $R$ etc.\ are used as in this theorem.

Theorem \ref{theo_genvertsum} provides a flexible tool to investigate posets and their down-set lattices systematically, because for every down-set $D \in \D(R)$, the poset $R$ is a g.v.s.\ of $R \vert_D$ and $R \setminus D$. With $D$ being the antichain of the minimal points of $R$, this approach has been used in \cite{aCampo_2019,aCampo_Erne_2019} for the enumeration of down-sets of posets and for the enumeration of posets with a certain characteristic.

We define for every $T \in \D(Q)$
\begin{align*}
\J_T(R) & \equivA \mysetdescr{ D \in \D(R) }{ D \cap Y = T }.
\end{align*}
Due to $Q = R \vert_Y$, the sets $\J_T(R), T \in \D(Q)$, form a partition of $\D(R)$. Therefore, 
\begin{align*}
h(R) & \eqA  \sum_{T \in \D(Q)} a_T(R), \\
\mytext{where} \quad
a_T(R) & \equivA \sum_{D \in \J_T(R)} \# \E_D(R) \quad \mytext{for all } T \in \D(Q).
\end{align*}
In the case of $X \subseteq \mydownarrow_R Y$, we have $\J_Y(R) = \{ X \cup Y \}$, hence
\begin{align} \label{aTR_DRXT}
a_Y(R) & \eqA \# \D ( R ).
\end{align}

\begin{definition} \label{def_PQR}
In what follows, $B^- \subseteq X$ is a fixed up-set of $P$ and $B^+ \subseteq Y$ is a fixed down-set of $Q$. We set $S^- \equiv P \vert_{B^-}$, $S^+ \equiv Q \vert_{B^+}$, and we assume that $\sigma : \D( S^+ ) \rightarrow \P( B^- )$ is a mapping with
\begin{align} \label{sigma_compatible}
\forall \; T \in \D(S^+) \mytext{: } & \sigma(T) \; \subseteq \; X \cap \mydownarrow_R T \; \subseteq \;  \mydownarrow_P \; \sigma(T).
\end{align}
\end{definition}
Because $B^-$ is an up-set of $P$, the poset $P$ is a g.v.s.\ of $P \setminus B^-$ and $S^-$; similarly, $Q$ is a g.v.s.\ of $S^+$ and $Q \setminus B^+$. For later use, we note that the first inclusion in \eqref{sigma_compatible} enforces 
\begin{align} \label{sigma_leer}
\sigma( \emptyset ) & \eqA \emptyset.
\end{align}

Firstly, we realize that the assumptions in Definition \ref{def_PQR} are not restrictive. For a given poset $R$, let $U$ be an up-set different from $R$ and $\emptyset$. We define $P \equiv R \setminus U$ and $Q \equiv R \vert_U$. We select for $B^+$ any non-empty down-set of $Q$ and $B^-$ as the up-set of $P$ created by the points of $P$ which are covered by points of $B^+$ in $R$. ($B^-$ can  be empty; we come back to this case at the end of the section.)  With the mapping
\begin{align*}
\sigma : \D(S^+) & \rightarrow \P(B^-), \\
T & \mapsto B^- \; \cap \; \mydownarrow_R T,
\end{align*}
all requirements in Definition \ref{def_PQR} are fulfilled. However, such a schematic choice of $S^+$, $S^-$, and $\sigma$ may be unfavorable. For the effective calculation of $h(R)$, they should be selected in such a way, that they match the structure of the recursion, as discussed at the beginning of Section \ref{sec_application}.

The isomorphisms in the following lemma are the key for the recursive approach; they are illustrated in the Figures \ref{fig_stepPyramid} and \ref{fig_HCC_DHCC} in Section \ref{sec_application}.

\begin{lemma} \label{lemma_beta_tau}
For every $T \in \D(S^+)$, the mapping
\begin{align} \label{tau_is_isom}
\tau_T : \J_T(R) & \rightarrow \bigcup_{\stackrel{U \in \D(S^-)}{\sigma(T) \subseteq U}} \J_U(P), \\ \nonumber
D & \mapsto D \setminus T
\end{align}
is an isomorphism with inverse $D' \mapsto D' \cup T$; in particular
\begin{align} \label{Jempty_isom}
\J_\emptyset(R) & \simeq \D(P).
\end{align}
Furthermore, for every $N \in \D(Q)$, the mapping
\begin{align} \label{beta_is_isom}
\beta_N : \mysetdescr{ D \in \D(R) }{ N \subseteq D } & \rightarrow \D( R \smdaR N ), \\ \nonumber
D & \mapsto D \smdaR N
\end{align}
is an isomorphism with inverse $D' \mapsto D' \cup \mydownarrow_R N$. In particular, for $T \in \D(Q)$ with $N \subseteq T$,
\begin{align} \label{isom_setminus}
\J_T(R) & \simeqA \J_{T \setminus N}( R \smdaR N )
\end{align}
\end{lemma}
\BP \eqref{tau_is_isom}: Let $T \in \D(S^+)$. Because $B^+$ is a down-set of $Q$, we have $\D(S^+) \subseteq \D(Q)$, and $\J_T(R)$ is well-defined. And because $P$ is a g.v.s.\ of $P \setminus B^-$ and $S^-$, also $\J_U(P)$ is well-defined for every $U \in \D(S^-)$.

Let $D \in \J_T(R)$. Because $P$ is the lower part of $R$, we have $D \setminus T \in \D(P)$, and due to the first inclusion in \eqref{sigma_compatible}, we even have $\sigma(T) \subseteq D \setminus T$, thus $\sigma(T) \subseteq ( D \setminus T ) \cap B^- $. The sets $\J_V(P)$, $V \in \D(S^-)$, form a partition of $\D(P)$; therefore, the set $D \setminus T$ belongs to the union on the right of \eqref{tau_is_isom}, and the mapping $\tau_T$ is well-defined.

Let $D'$ belong to the union on the right of \eqref{tau_is_isom}. By case discrimination, we show that $D \equiv D' \cup T$ is a down-set of $R$. Let $x \in D$ and $y \in R$ with $y \leq_R x$:
\begin{itemize}
\item $x \in D', y \in P$: $D' \in D(P)$ yields $y \in D$.
\item $x \in D', y \in Q$: impossible due to $\leq_R \cap \; ( Y \times X) = \emptyset$.
\item $x \in T, y \in P$: The second inclusion in \eqref{sigma_compatible} delivers 
\begin{equation*}
y \; \inA \mydownarrow_P \; \sigma(T) \; \; \subseteq \; \; D'.
\end{equation*}
\item $x \in T, y \in Q$: $T \in \D(S^+) \subseteq \D(Q)$ yields $y \in D$.
\end{itemize}
Therefore, $D$ is a down-set of $R$, and $D \in \J_T(R)$ follows. We conclude that $\tau_T$ has the inverse $D' \mapsto D' \cup T$. Isomorphism follows, because $\tau_T$ and its inverse are both homomorphisms with respect to ``$\subseteq$''.

\eqref{Jempty_isom}: Follows with \eqref{tau_is_isom} and \eqref{sigma_leer}, because the sets $\J_U(P)$, $U \in \D(S^-)$, form a partition of $\D(P)$.

\eqref{beta_is_isom}: It is easily seen that $D \smdaR M$ is indeed a down-set of $R \smdaR M$ for every $D \in \D(R)$ and every subset $M \subseteq X \cup Y$.

Let $N \in \D(Q)$. We show that the indicated inverse of $\beta_N$ is well-defined. Let $P' \equiv P \smdaR N$, $Q' \equiv Q \setminus N$, $R' \equiv R \smdaR N$. The set $X \smdaR N$ is a down-set of $R'$ and $Y \setminus N$ is an up-set of $R'$; the poset $R'$ is thus a g.v.s.\ of $P'$ and $Q'$. Let $D \in \D( R' )$. According to \eqref{formel_downsets}, there exists a unique $T' \in \D(Q')$ and a unique $D' \in \D( P' \smda{R'} T' )$ with $D = D' \cup \mydownarrow_{R'} T'$. We have $N \cup T \in \D(Q)$ and additionally $P' \smda{R'} T' = P \smdaR (N \cup T)$, and \eqref{formel_downsets} delivers that indeed $D \cup \mydownarrow_R N = D' \cup \mydownarrow_R (N \cup T)$ is a down-set of $R$ containing $N$.

Isomorphism follows again because $\beta_N$ and its inverse are both homomorphisms with respect to ``$\subseteq$''. \eqref{isom_setminus} is a direct consequence, because of $T \setminus N \in \D( Q' )$.

\EP

From the definitions, it is directly understandable that the coefficient $a_T(R)$ is not affected by $Q \setminus T$. The proof of this useful observation is technical:
\begin{lemma} \label{lemma_aTRStr_aTR}
Let $Y' \in \D(Q)$. With $Q' \equiv Q \vert_{Y'}$, $R' \equiv R \vert_{X \cup Y'}$, we have $a_T(R') = a_T(R)$ for every $T \in \D(Q')$.
\end{lemma}
\BP Because of $Y' \in \D(Q)$ and $X \cup Y' \in \D(R)$, we have $\D(Q') \subseteq \D(Q)$ and $\D(R') \subseteq \D(R)$. Moreover, $\D(Q')$ is a down-set of $\D(Q)$, and $R'$ is a g.v.s.\ of $P$ and $Q'$. 

Let $T \in \D(Q')$. The equation $a_T(R') = a_T(R)$ will be a direct consequence of
\begin{align} \label{JTRStr_JTR}
\J_T(R') & = \J_T(R), \\ \label{EDRStr_EDR}
\mytext{and } \quad \E_D(R') & = \E_D(R) \quad \mytext{for all } D \in \J_T(R').
\end{align}
Because of $T \subseteq Y'$, we have $\leq_{R'} \capA ( X \times T )\eqA \leq_R \capA ( X \times T )$, hence $X \cap \mydownarrow_{R'} T = X \cap \mydownarrow_R T$. Because additionally $Y' \cap \mydownarrow_{R'} T = T = Y \cap \mydownarrow_R T$, 
\begin{align} \label{eq_RStirchT_RT}
\mydownarrow_{R'} T & \eqA \mydownarrow_R T,
\end{align}
and applying \eqref{formel_downsets} twice yields \eqref{JTRStr_JTR}:
\begin{align*}
\J_T(R') & \eqA \mysetdescr{ D \; \cup \downarrow_{R'} T }{ D \in \D( P \smda{R'} T) } \\
& \eqA \mysetdescr{ D \; \cup \downarrow_R T }{ D \in \D( P \smdaR T) } \eqA \J_T(R).
\end{align*}
Now let $D \in \D(Q')$. $\mydownarrow_{\D(Q)} D$ is a down-set of $\D(Q)$ with $D \in \D(Q')$, and because $\D(Q')$ is a down-set of $\D(Q)$, we have $\mydownarrow_{\D(Q)} D \subseteq \D(Q') $, hence $\mydownarrow_{\D(Q)} D \subseteq \mydownarrow_{\D(Q')} D$. We conclude $\mydownarrow_{\D(Q)} D = \mydownarrow_{\D(Q')} D$, because $Q'$ is a sub-poset of $Q$. All togehter, application of \eqref{eq_RStirchT_RT} yields
\begin{align*}
\mydownarrow_{R'} O & \eqA \mydownarrow_{R} O \quad \mytext{for all } O \in \mydownarrow_{\D(Q')} D \eqA  \mydownarrow_{\D(Q)} D,
\end{align*}
implying $\mydownarrow_{R'} D \eqA \mydownarrow_{R} D$. Now \eqref{EDRStr_EDR} results, because due to \eqref{formel_downsets},
\begin{align*}
\E_D(R') & \eqA \bigcup_{O \in \mydownarrow_{\D(Q')} D} \mysetdescr{ F \; \cup \downarrow_{R'} O }{ F \in \D( P \smda{R'} O ), \; F \subseteq X \cap \mydownarrow_{R'} D }, \\
\E_D(R) & \eqA \bigcup_{O \in \mydownarrow_{\D(Q)} D} \mysetdescr{ F \; \cup \downarrow_{R} O }{ F \in \D( P \smda{R} O ), \; F \subseteq X \cap \mydownarrow_{R} D }.
\end{align*}

\EP

In the following theorem, it is described how the coefficients $a_T(R)$ with $T \in \D(S^+)$ can be determined recursively:
\begin{theorem} \label{theo_formel_a}
Let $M$ be the set of the minimal points of $Q$ and let $P^+ \equiv R \vert_{X \cup B^+}$. Then, for all $T \in \D(S^+)$
\begin{align} \label{formel_a}
a_T(R) & \; = \; \sum_{\stackrel{U \in \D(S^-)}{\sigma(T) \subseteq U}} a_U(P)
\; \;  + \; \; \sum_{\mu = 1}^{\# ( M \cap T)} (-1)^{\mu-1} \sum_{\stackrel{N \subseteq M \cap T}{\# N = \mu}} a_{ T \setminus N }(P^+ \smdaPp N).
\end{align}
In particular,
\begin{align} \label{a_empty_R}
a_\emptyset(R) & = h( P ).
\end{align}
\end{theorem}
\BP We prove the following two equations which immediately yield \eqref{formel_a}:
\begin{align} \label{summe_ERDempty}
\sum_{D \in \J_T(R)} \# \left( \E_D(R) \cap \J_\emptyset(R) \right) & = 
\sum_{\stackrel{U \in \D(S^-)}{\sigma(T) \subseteq U}} a_P(U), \\ \label{summe_ERDnotempty}
\sum_{D \in \J_T(R)} \# \left( \E_D(R) \setminus \J_\emptyset(R) \right) & = \sum_{\mu = 1}^{\# ( M \cap T)} (-1)^{\mu-1} \sum_{\stackrel{N \subseteq M \cap T}{\# N = \mu}} a_{ T \setminus N }(P^+ \smdaPp N).
\end{align}

\eqref{summe_ERDempty}: We start with the proof of $\E_D(R) \cap \J_\emptyset(R) = \E_{\tau_T(D)}(P)$. Let $E \in \E_D(R) \cap \J_\emptyset(R)$. $E \in \D(R)$ and $E \cap Y = \emptyset$ means $E \in \D(P)$, and $E \subseteq D$ additionally yields $E = E \setminus T \subseteq D \setminus T= \tau_T(D)$. According to Lemma \ref{lemma_beta_tau}, $\tau_T(D)$ is a down-set of $P$, and $E \in \E_{\tau_T(D)}(P)$ follows. On the other hand, $E' \in \E_{\tau_T(D)}(P)$ yields $E' \in \D(R)$ with $E' \subseteq \tau_T(D) = D \setminus T \subseteq D$ and $E' \cap Y = \emptyset$, hence $E' \in \E_D(R) \cap \J_\emptyset(R)$.

Now Lemma \ref{lemma_beta_tau} yields
\begin{align*}
& \sum_{D \in \J_T(R)} \# \left( \E_D(R) \cap \J_\emptyset(R) \right) =
\sum_{D \in \J_T(R)} \# \E_{\tau_T(D)}(P) \\
\stackrel{\eqref{tau_is_isom}}{=} &
\sum_{\stackrel{U \in \D(S^-)}{\sigma(T) \subseteq U}} 
\sum_{D' \in \J_U(P)} \# \E_{D'}(P) 
= \sum_{\stackrel{U \in \D(S^-)}{\sigma(T) \subseteq U}} a_U(P).
\end{align*}

\eqref{summe_ERDnotempty}: For $T = \emptyset$, the right side of the formula is zero, and also the left side is zero due to $D \in \J_\emptyset(R) \Rightarrow \E_D(R) \subseteq \J_\emptyset(R)$. 

For $T \not= \emptyset$, we have $M \cap T \not= \emptyset$. Let $\emptyset \not= N \subseteq M \cap T$. Then, for $D \in \J_T(R)$,
\begin{align*}
\E_D^N(R) & \equivA 
\mysetdescr{ E \in \E_D(R) }{ N \subseteq E } \\
& \eqA \mysetdescr{ E \in \D(R) }{ N \subseteq E \subseteq D} \\
& \stackrel{\eqref{beta_is_isom}}{\simeqA}
\mysetdescr{ E \in \D( R \smdaR N ) }{ E \subseteq D \smdaR N} \\
& \eqA \E_{D \smdaR N }( R \smdaR N ).
\end{align*}
Now \eqref{isom_setminus} yields with $R' \equiv R \smdaR N$,
\begin{align*}
\sum_{D \in \J_T(R)} \# \E_D^N(R) & = 
\sum_{D' \in \J_{T \setminus N}(R')} \# \E_{D'}( R' ) = 
 a_{T \setminus N}( R' ).
\end{align*}

For every $D \in \J_T(R)$, the set $\E_D(R) \setminus \J_\emptyset(R)$ contains exactly the sets $E \in \E_D(R)$ with $m \in E$ for an $m \in M \cap T$. Therefore,
\begin{align*}
& \sum_{D \in \J_T(R)} \# \left( \E_D(R) \setminus \J_\emptyset(R) \right) \\
= & \sum_{D \in \J_T(R)} \sum_{\mu = 1}^{\# ( M \cap T)} (-1)^{\mu-1} \sum_{\stackrel{N \subseteq M \cap T}{\# N = \mu}} \# \E_D^N(R) \\
\\
= & \sum_{\mu = 1}^{\# ( M \cap T)} (-1)^{\mu-1} \sum_{\stackrel{N \subseteq M \cap T}{\# N = \mu}} a_{T \setminus N}( R' ),
\end{align*}
and \eqref{summe_ERDnotempty} follows, because Lemma \ref{lemma_aTRStr_aTR} delivers $a_{T \setminus N}( R' ) = a_{T \setminus N}( P^+ \smdaPp N)$ for $T \in \D(S^+)$.

\eqref{a_empty_R}: For $T = \emptyset$, the double sum on the right of \eqref{formel_a} is zero, hence
\begin{align*}
a_\emptyset(R) & 
\stackrel{\eqref{formel_a} \eqref{sigma_leer}}{=}
\sum_{U \in \D(S^-)} a_U(P) \; = \; h(P).
\end{align*}
because $P$ is a g.v.s.\ of $P \setminus B^-$ and $S^-$.

\EP

In Section \ref{sec_application}, we frequently work with structures as in the following corollary:

\begin{corollary} \label{coro_formel_h}
Assume
\begin{align*}
& Q \mytext{has a minimum point} \bot, \\
& \sigma : \D(S^+) \rightarrow \D(S^-) \mytext{is an isomorphism}.
\end{align*}
Then, setting $P^+ \equiv R \vert_{X \cup B^+}$ again, 
\begin{align} \label{formel_h}
\begin{split}
h(R) \; \;  = \; \; h(  P^+ \smdaPp \bot )
\; \; & + \; \;
\sum_{D \in \D(Q) \setminus \D(S^+)} a_D(R) \\
& + \; \; \sum_{T \in \D(S^-)} \left( \# \mydownarrow_{\D(S^-)} T \right) \cdot a_T(P).
\end{split}
\end{align}
\end{corollary}
\BP Observing that the double-sum on the right of \eqref{formel_a} is zero for $T = \emptyset$, we get
\begin{align*}
& \sum_{T \in \D(S^+)} a_T(R) \\
= & \sum_{T \in \D(S^+)} \left( \sum_{\stackrel{U \in \D(S^-)}{\sigma(T) \subseteq U}} a_U(P)
+ \sum_{\emptyset \not= N \subseteq \{ \bot \} \cap T} a_{T \setminus N}( P^+ \smdaPp N ) \right) \\
= & \sum_{T \in \D(S^-)} \left( \# \mydownarrow_{\D(S^-)} T \right) \cdot a_T(P)
\; \; + \sum_{T \in \D(S^+) \setminus \{ \emptyset \}} a_{ T \setminus \{ \bot \} }(P^+ \smdaPp \bot).
\end{align*}
Because $P^+ \smdaPp \bot$ is a g.v.s.\ of $P \smdaPp \bot $ and $S^+ \setminus \{ \bot \}$, the right sum is
\begin{align*}
\sum_{T \in \D(S^+) \setminus \{ \emptyset \}} a_{ T \setminus \{ \bot \} }(P^+ \smdaPp \bot)
\; = \; & \sum_{T \in \D(S^+ \setminus \{ \bot \} ) } a_{T}( P^+ \smdaPp \bot ) \\
\; = \; & h( P^+ \smdaPp \bot ).
\end{align*}

\EP

However, Formula \eqref{formel_h} is of limited value, because for the calculation of the $a$-coefficients, we still need Formula \eqref{formel_a} from Theorem \ref{theo_formel_a}.

We have mentioned after Definition \ref{def_PQR} that we can run the recursive approach always with $B^+$ being any non-empty down-set of $Q$ and $B^-$ being the up-set of $P$ created by the points of $P$ covered by points of $B^+$ in $R$. For these choices, $B^-$ may be empty. (An example is the poset $R$ in Figure \ref{fig_iii_i}; take the two large Lambdas as $P$ and the remaining part as $Q$, and define $B^+$ as the singleton containing the minimum point of $Q$ only.) In this case, $S^-$ is the empty poset, $\sigma(T) = \emptyset$ for all $T \in \D(S^+) = \P(B^+)$, and the isomorphism $\tau_T$ in \eqref{tau_is_isom} reduces to
\begin{align*}
\tau_T : \J_T(R) & \rightarrow \J_\emptyset(P) 
\stackrel{\eqref{Jempty_isom}}{\eqA}
\D(P), \\
D & \mapsto D \setminus T.
\end{align*}
Indeed, if no point of $P$ is covered by a point of $B^+$, then $R \vert_{X \cup B^+} = P + S^+$, hence $\J_T(R) \simeq \D(P) \times \{ T \} \simeq \D(P)$ for every $T \in \D(S^+) = \P(B^+)$.

\section{Application} \label{sec_application}

According to Theorem \ref{theo_formel_a}, we need for the calculation of the coefficients $a_T(R)$ with $T \in \D(S^+)$ the values of $a_U(P)$ with $U \in \D(S^-)$ and the coefficients $a_{ T \setminus N }(P^+ \smdaPp N)$ for all non-empty sets $N$ of minimal points of $S^+$. Additionally, we need the coefficients $a_D(R)$ with $D \in \D(Q) \setminus \D(S^+)$ for the final calculation of $h(R)$. For a given poset $R$, the choice of $P$, $Q$, and $S^+$ shifts the balance between these three types of calculation:
\begin{itemize}
\item A large sub-poset $S^+$ of $Q$ reduces the number of down-sets $D \in \D(Q) \setminus \D(S^+)$ for which $a_D(R)$ has to be calculated separately (even to zero for $S^+ = Q$), but makes the calculation of the coefficients $a_{ T \setminus N }(P^+ \smdaPp N)$ in the second sum in \eqref{formel_a} more demanding.
\item A small sub-poset $P$ of $R$ makes the determination of the coefficients $a_U(P)$ with $U \in \D(S^-)$ easier in the first sum in \eqref{formel_a}, but puts a larger burden on at least one of the two other types of calculation.
\item The number of minimal points of $S^+$ affects exponentially the number of terms in the second sum in \eqref{formel_a}.
\end{itemize}
Our approach fits thus best to posets $R$ with the following properties:
\begin{itemize}
\item $P$ and $P^+$ have a simple structure or are closely related to $R$;
\item $a_D(R)$ can easily be calculated for all $D \in \D(Q) \setminus \D( S^+ )$;
\item the number of minimal points of $S^+$ is small.
\end{itemize}
For such posets, the calculation of $h( R )$ is possible with ordinary table calculation, as we will see in this section. We work with posets $R = W \times C_k$ for different posets $W$ in Section \ref{subsec_LambdaDiamond} (including $W = C_3 \times C_3$), and with $R = \H( C_2, C_k ) \simeq \H( C_{k-1}, C_3 )$ in Section \ref{subsec_HC2Ck}. The posets $S^+$ and $S^-$ are always isomorphic, and the mapping $\sigma$ fulfilling \eqref{sigma_compatible} is induced by the respective isomorphism $\iota : S^+ \rightarrow S^-$:
\begin{align*}
\forall \; T \in \D(S^+) \mytext{: } \sigma(T) & \equivA \iota[ T ].
\end{align*}
In Section \ref{subsec_LambdaDiamond}, we have $S^+ = Q$ in all cases, but in Section \ref{subsec_HC2Ck}, we will see that a different choice of $S^+$ can even be beneficial.

For almost all posets $R$ in this section, we have $X \subseteq \mydownarrow_R Y$. (The exception is the poset in the lower part of Figure \ref{fig_iii_i}). For these, equation  \eqref{aTR_DRXT} delivers $a_Q(R) = \# \D(R) = \# \H(R,C_2)$, and as a by-product, we get the number of {\em surjective} homomorphisms $R \rightarrow C_3$ via
\begin{align} \label{anz_surjHoms}
h(R) \; - \; 3 \cdot a_Q(R) \; + \; 3.
\end{align}

Due to the general isomorphism $\H(P_1, \H(P_2,P_3) ) \simeq \H( P_1 \times P_2, P_3)$ and $\H(C_2,C_2) \simeq C_3$, we have $h(R) = \# \D( R \times C_2)$ for every poset $R$. In particular, the recursive approach gives raise to new ways to calculate Dedekind  numbers, because the $(k+1)$\textsuperscript{th} Dedekind number is $\# \D(C_2^{k+1} ) = h( C_2^k )$.

\subsection{$ R = \WCk$} \label{subsec_LambdaDiamond}

For every poset $W$, the product $\WCk$ fits into the frame of Definition \ref{def_PQR} via
\begin{align*}
P & \equivA \WCkm \\
Q & \equivA W \times \{k\}, \\
S^- & \equivA W \times \{k-1\}, \\
S^+ & \equivA W \times \{k\}, \\
\iota(w,k) & \equivA (w,k-1) \\
P^+ & \eqA W \times C_k.
\end{align*}
In order to avoid unnecessary formalism, we identify $S^-$ and $S^+$ with $W$. Due to \eqref{aTR_DRXT}, we always have $a_W(W \times C_k) = \# \D( W \times C_k )$.

In the following sections, we determine $h( C_n \times C_k)$, $h( \LCk )$, $h( \LOCk )$, and $h( \CCdreik )$, where $\Lambda \equiv A_2 \oplus A_1$ is the poset with $\Lambda$-shaped diagram and $\lozenge \equiv C_2 \times C_2$ is the diamond. As a by-product of the calculation of $h( C_3 \times C_3 \times C_k )$, we get $h( W \times C_k)$ for eight additional posets $W$.

\subsubsection{$ \CCnk $} \label{subsubsec_Cn}


In order to simplify notation, we write $a_j(\CCnk)$ instead of $a_{\myNk{j}}( \CCnk )$ for all $j \in \myNk{n}_0$. For $k=1$, $a_j(C_n) = j+1$ and $ h(C_n) = \frac{(n+1)(n+2)}{2}$ are trivial.

\begin{theorem}
For $k \geq 2$,
\begin{align*}
a_0(\CCnk) & \eqA h(\CC{n}{k-1}), \\ \nonumber
\forall j \in \Nn \mytext{: } a_j(\CCnk) & \eqA a_{j-1}( \CC{n-1}{k} )
+ \sum_{i=j}^{n} a_{i}(\CC{n}{k-1}), \\
\mytext{hence} \quad h( \CCnk ) & \eqA h( \CC{n-1}{k} ) + \sum_{i=0}^{n} ( i+1 ) \cdot a_{i}(\CC{n}{k-1}).
\end{align*}
\end{theorem}
\BP The first equation is due to \eqref{Jempty_isom}. Let $j \in \myNk{n}$. The symbols in the left sum of  \eqref{formel_a} become
\begin{align*}
T & \eqA \myNk{j}, \\
\mysetdescr{ U \in \D(S^-) }{ \iota[T] \subseteq U } & = \mysetdescr{ \myNk{i} }{ j \leq i \leq n }.
\end{align*}
The right sum reduces to the single term $a_{T \setminus \{ \bot \}}( P^+ \smdaPp \bot )$ with
\begin{align*}
\bot & \eqA 1, \\
T \setminus \{ \bot \} & \simeqA C_{j-1}, \\
P^+ \smdaPp \bot & \simeqA C_{n-1} \times C_k,
\end{align*}
and the formula for $a_j(\CCnk)$ with $j \in \myNk{n}$ follows. The last formula is \eqref{formel_h}.

\EP

The coefficients $h( \CCnk )$ and $a_j(\CCnk)$ are shown in Figure \ref{fig_Tabelle_CnCk} in the Appendix for $n \in \myNk{4}$ and $k \in \myNk{5}$.

\subsubsection{$\Lambda \times C_k$ and $\lozenge \times C_k$} \label{subsubsec_LambdaRaute}

\begin{figure} 
\begin{center}
\includegraphics[trim = 70 630 370 70, clip]{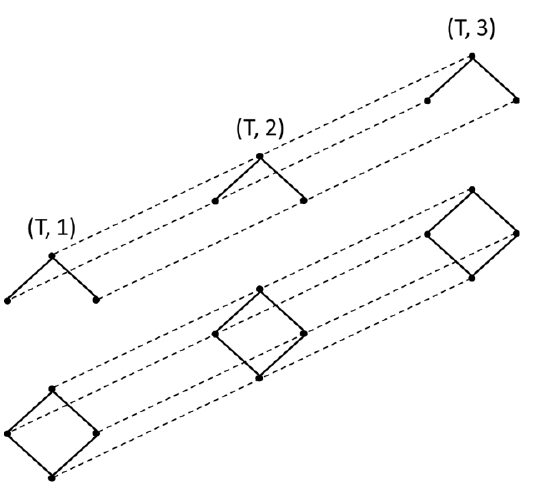}
\caption{\label{fig_WxC3} The posets $\Lambda \times C_3$ and $\lozenge \times C_3$. }
\end{center}
\end{figure}

The posets $\Lambda \times C_3$ and $\lozenge \times C_3$ are shown in Figure \ref{fig_WxC3}. We start with $\Lambda$. Even if one of our results is a closed formula for $h( \LCk )$, we remain interested in the coefficients $a_T(\LCk )$ and derive formulas for them, because we need them for $\lozenge$ and $C_3 \times C_3$.

\begin{figure} 
\begin{center}
\includegraphics[trim = 70 450 230 80, clip]{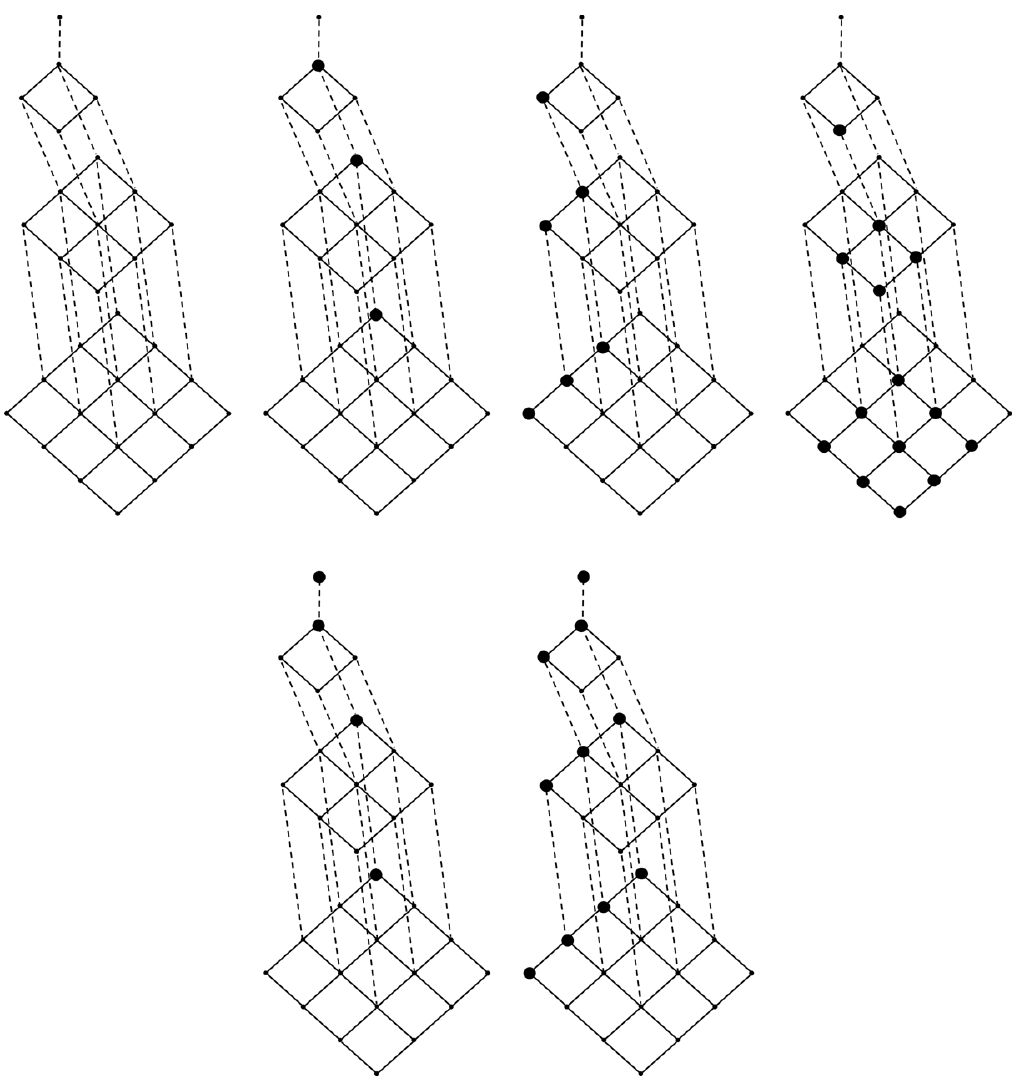}
\caption{\label{fig_stepPyramid} Top row, from left to right: $\D(\LC{3})$, $\J_\Mellerr( \LC{3})$, $\J_\Mell( \LC{3})$, and $\J_\emptyset( \LC{3}) \simeq \D(\LC{2})$. Bottom row: illustration of the isomorphisms $\beta_\Mellerr$ and $\beta_\Mell$ from Lemma \ref{lemma_beta_tau}; explanations in text.}
\end{center}
\end{figure}

We denote with $\top$ the maximum point of $\Lambda$ and with $\ell, r$ the minimal points. In $\D(\LCk)$, the set of down-sets $D$ with $( \top, j ) \in D$ and $( \top, j+1 ) \notin D$ is isomorphic to $C_{k+1-j} \times C_{k+1-j}$. All together, $\D(\LCk)$ looks like a step pyramid, as shown in Figure \ref{fig_stepPyramid} for $\LC{3}$. We denote with $f$ the ``floor'' of the pyramid, starting with $f = 0$ for the ground floor and ending with $f = k$ for the ``antenna'' $\LCk$ on top.

The set $\J_\Lambda(\LCk)$ contains the antenna only, and $\J_\Mellerr(\LCk)$ consists of the points of the ``backbone'' of the pyramid marked in the second drawing in the top row of Figure \ref{fig_stepPyramid}:
\begin{align*}
\J_\Mellerr( \LCk ) & \eqA \mysetdescr{ D_f }{ f \in \Nke_0 }, \\
\mytext{with} \quad D_f & \equivA ( \LCk ) \setminus \mysetdescr{ (\top,j) }{ f+1 \leq j \leq k }
\end{align*}
$\J_\Mell(\LCk)$ and $\J_\Merr(\LCk)$ are the respective back sides of the pyramid without the backbone, and $\J_\emptyset(\LCk) \simeq \D( \LCkm)$ contains the inner points and the front points of the pyramid, as shown in Figure \ref{fig_stepPyramid}.

The isomorphisms \eqref{tau_is_isom} in Lemma \ref{lemma_beta_tau} become
\begin{align*}
\J_\Mellerr(\LCk) \simeqA & \; \; \; \; \; \J_\Mellerr( \LCkm ) \cup \J_{\Lambda}(\LCkm ) \\
\J_\Mell(\LCk) \simeqA & \; \; \; \; \; \J_\Mell( \LCkm ) \\ & \cupA \J_\Mellerr( \LCkm ) \cup \J_{\Lambda}(\LCkm ),
\end{align*}
and indeed, the backbone plus the antenna of $\D( \LCkm )$ and the backbone of $\D( \LCk )$ are both isomorphic to $C_k$, and adding the back side $\J_\Mell( \LCkm )$ of the pyramid $\D( \LCkm )$ yields a lattice isomorphic to $\J_\Mell( \LCk )$ (cf.\ the upper part of Figure \ref{fig_stepPyramid}).  Finally, due to the isomorphisms $\beta_\Mellerr$ and $\beta_\Mell$ in Lemma \ref{lemma_beta_tau},
\begin{align*}
\mysetdescr{ D \in \D( \LCk ) }{ (\ell,k), (r,k) \in D } & \simeqA \D( C_k ), \\
\mysetdescr{ D \in \D( \LCk ) }{ (\ell,k) \in D } & \simeqA \D( C_2 \times C_k ),
\end{align*}
as shown in the lower part of Figure \ref{fig_stepPyramid}.

The number of points contained in the $f$\textsuperscript{th} floor is $(k+1-f)^2$, hence
\begin{align*}
a_{\Lambda}( \LCk ) & \eqA \sum_{f=0}^k ( k+1-f )^2 \eqA \frac{ (k+1)(k+2)(2k+3)}{6}.
\end{align*}

For the calculation of $a_\Mellerr( \LCk )$, we realize that $\E_{D_f}(\LCk)$ consists of the points belonging to the floors $0 \ldots f$. Each floor contributes thus $(k-f)$-times to $a_\Mellerr( \LCk )$, hence
\begin{align*}
a_\Mellerr( \LCk ) & \eqA \sum_{f=0}^{k-1} ( k - f ) (k + 1 - f )^2 \eqA \frac{k (k+1)(k+2)(3k+5)}{12}.
\end{align*}

Finally, it is easily seen that the contribution of the $f$\textsuperscript{th} floor to the number $a_\emptyset(\LCk)$ is 
\begin{align*}
F_f(k) & \equivA \sum_{i=1}^{k-f} \sum_{j=1}^{k-f} \sum_{\phi=0}^f (i+\phi)(j+\phi) \\
& \eqA \quad \Big( f (f+2) + 3 (k+2)^2 \Big) \cdot \frac{(f+1)(k+1-f)^2}{12},
\end{align*}
and summing up $F_f(k)$ from $0$ to $k$ delivers $a_\emptyset(\LCk)$. In particular, due to \eqref{a_empty_R},
\begin{align*}
h(\LCk) & \eqA \sum_{f=0}^{k+1} F_f(k+1).
\end{align*}
The missing numbers $a_\Mell(\LCk) = a_\Merr(\LCk)$ are given by
\begin{align*}
\frac{1}{2} \cdot \Big( h( \LCk ) - a_{\Lambda}( \LCk ) - a_\Mellerr( \LCk ) - a_\emptyset( \LCk ) \Big).
\end{align*}

We come to $\lozenge = C_2 \times C_2$. We denote with $\top$ and $\bot$ the maximum and minimum point of $\lozenge$ and with $\ell, r$ the remaining points. Due to 
\begin{align*}
\left( \LOCk \right) \smda{\LOCk} \bot & \simeqA \LCk,
\end{align*}
Theorem \ref{theo_formel_a} yields
\begin{align*}
a_\emptyset(\LOCk) & \eqA \quad h( \LOCkm ), \\
a_{\{ \bot \}}(\LOCk) & \eqA \quad h( \LOCkm ) - a_\emptyset(\LOCkm) + h( \LCk ), \\
a_\Mellb(\LOCk) & \eqA \quad a_\Mellb(\LOCkm) + a_\Mellerrb(\LOCkm)  \\
& \quad \; \; + a_\lozenge(\LOCkm) + a_\Mell(\LCk), \\
a_\Merrb(\LOCk) & \eqA \quad a_\Mellb(\LOCk), \\
a_\Mellerrb(\LOCk) & \eqA \quad a_\Mellerrb(\LOCkm) + a_\lozenge(\LOCkm) \\
& \quad \; \; + a_\Mellerr(\LCk), \\
a_\lozenge( \LOCk ) & \eqA \quad a_\lozenge(\LOCkm) + a_\Lambda(\LCk).
\end{align*}

For $k \in \myNk{5}$, the values of $h( \LCk ), a_T( \LCk ), h( \LOCk )$, and $a_T( \LOCk )$ are shown in Figure \ref{fig_LambdaLozenge} in the Appendix. We have $a_\Lambda(\LCk) = \# \D( \LCk )$ and $a_\lozenge(\LOCk) = \# \D( \LOCk )$.

\subsubsection{$ C_3 \times C_3 \times C_k$ } \label{subsubsec_C3C3}

\begin{figure}
\begin{center}
\includegraphics[trim = 70 410 230 70, clip]{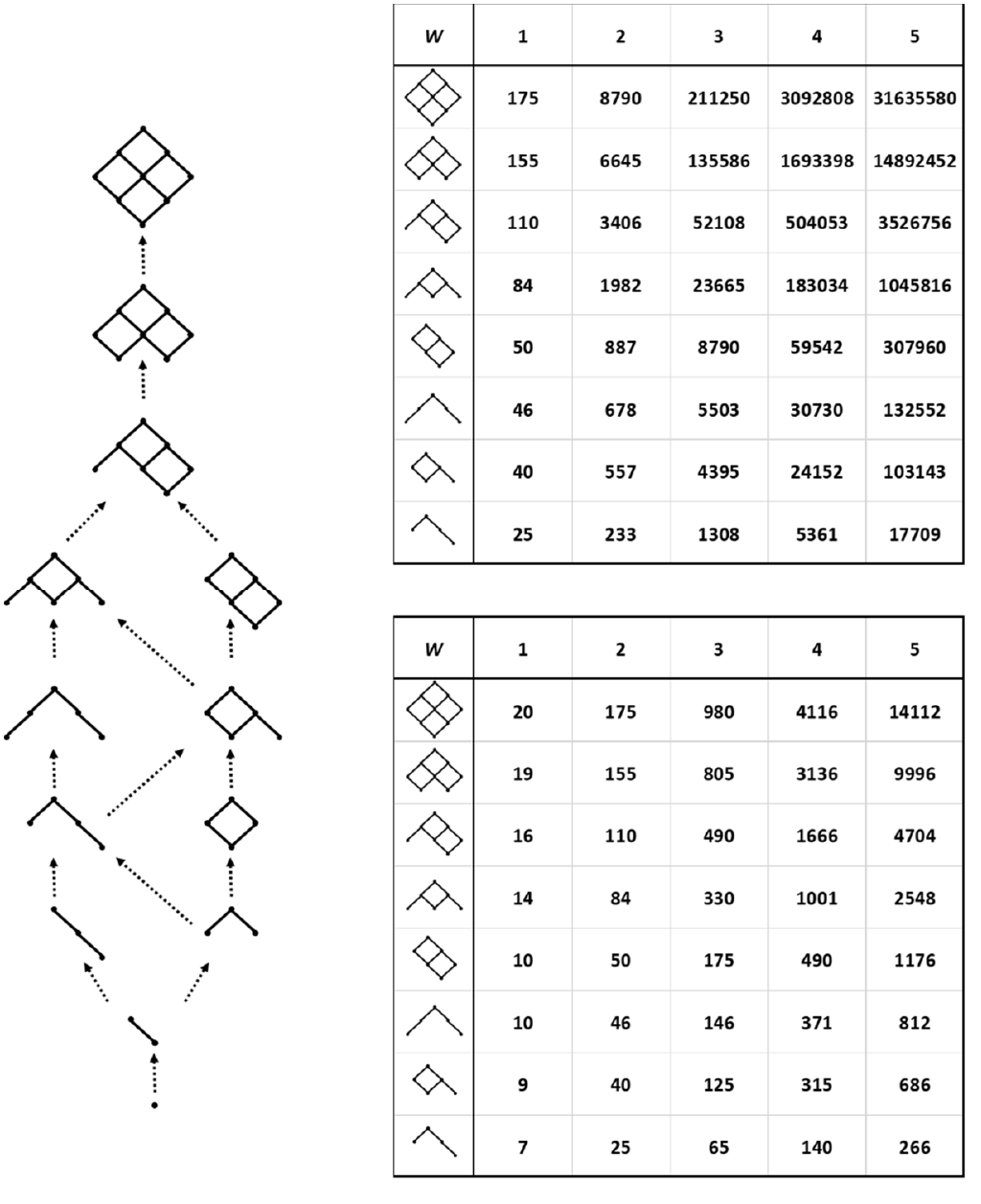}
\caption{\label{fig_C3C3C3C3} The posets $W$ required for the calculation of $h(\CCdreik)$, and the numbers $h(W \times C_k)$ (top table) and $\# \D( W \times C_k)$ (bottom table) for $k \in \myNk{5}$. For the chains, $\Lambda$, and $\lozenge$, see Figures \ref{fig_Tabelle_CnCk} and \ref{fig_LambdaLozenge} in the Appendix.}
\end{center}
\end{figure}

Due to the direct calculation of the coefficients $a_T( \LCk )$, only a single recursive step was required in the calculation of the coefficients $a_T( \LOCk )$. The case $\CCdreik$ is more demanding. We have to step recursively through the posets in the diagram in Figure \ref{fig_C3C3C3C3}. An arrow upwards from a poset $V$ to a poset $W$ in this transitive diagram indicates that the $a$-coefficients of $V \times C_k$ are required for the calculation of the $a$-coefficients of $W \times C_k$. In the figure, also the values of $h(W \times C_k)$ and $\# \D(W \times C_k )$ are shown for $k \in \myNk{5}$. For the chains and $\Lambda$, $\lozenge$, see Figures \ref{fig_Tabelle_CnCk} and \ref{fig_LambdaLozenge} in the Appendix.

\begin{figure} 
\begin{center}
\includegraphics[trim = 70 620 365 70, clip]{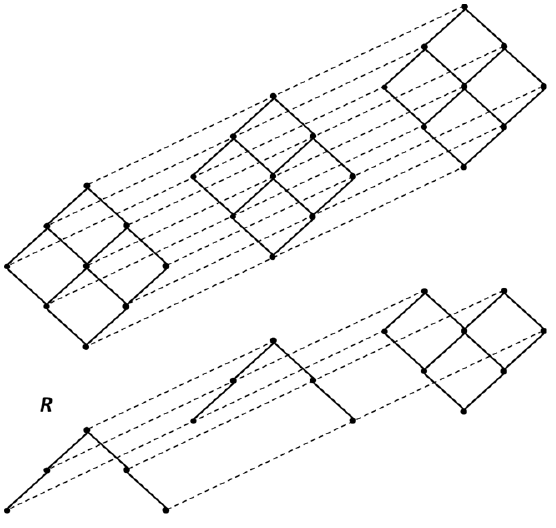}
\caption{\label{fig_iii_i} The poset $C_3^3$ and the poset $R$ required for the enumeration of the homomorphisms $\xi \in \H( C_3^3, C_3 )$ with $\xi(i, \ldots , i) = i$ for all $i \in \myNk{3}$. }
\end{center}
\end{figure}

At the end of the introduction, we mentioned machine learning and ensemble based systems as application of monotone ternary functions. Here, the {\em surjective} homomorphisms from $C_3^k$ to $C_3$ are of particular interest. We calculate their number with formula \eqref{anz_surjHoms}. Additional regularity is introduced by demanding that a summary rule $s$ has to respect an unanimous decision of the experts, i.e., $s(i, \ldots , i) = i$ for all $i \in \myNk{3}$. (Because we work with monotone summary rules only, this postulate is equivalent to ``$\min r \leq s(r) \leq \max r$ for all $r \in C_3^k$''.) The numbers of these homomorphisms are shown in the following table for $k = 1, 2, 3$.
\newline

\begin{tabular}{| c | c | c | c |}
\hline
$ k $ & $h( C_3^k )$ & surjective & $\xi(i, \ldots ,i) = i$ \\
\hline \hline
1 & 10 & 1 & 1 \\
2 & 175 & 118 & 64 \\
3 & 211250 & 208313 & 116211 \\
\hline
\end{tabular}
\newline
\newline

For $k = 3$, the number of the homomorphisms $\xi \in \H( C_3^3, C_3)$ with $\xi(i, i , i) = i$ for all $i \in \myNk{3}$ has been calculated as follows. With $R$ being the poset shown in the lower part of Figure \ref{fig_iii_i}, $h(R) = 46540$ is the number of homomorphisms $\xi \in \H( C_3^3, C_3)$ with $\xi(1,1,1) =  \xi(2,2,2) = 1$ and $\xi(3,3,3) = 3$. (The number has been calculated by applying the recursive method on the dual of $R$.) From these, $489 = \# \D(R) = \# \H(R,C_2)$ (calculated via \eqref{formel_downsets}) have $2$ not in their image; therefore, we have $46051$ surjective homomorphisms from $C_3^3$ to $C_3$ with $\xi(1,1,1) =  \xi(2,2,2) = 1$ and $\xi(3,3,3) = 3$. This is also the number of surjective homomorphisms from $C_3^3$ to $C_3$ with $\xi(1,1,1) = 1$ and $\xi(2,2,2) = \xi(3,3,3) = 3$, and $116211$ results as number of homomorphisms $\xi \in \H( C_3^3, C_3)$ with $\xi(i, i , i) = i$ for all $i \in \myNk{3}$.


\subsection{$R = \H( C_2, C_k ) \simeq \H(C_{k-1}, C_3)$} \label{subsec_HC2Ck}

As pointed out in the introduction, $h( C_3^k )$ is the number of the different ways how ternary rankings of an object by $k$ experts can be summarized by monotone summary rules in ensemble based systems. For several important summary rules like the different versions of {\em majority voting} \cite{Polikar_2006}, the resulting summary does not depend on who of the experts gave which rank; here, the summary rule has to be a {\em symmetric function}. In this case, we can order the rankings of the experts in non-decreasing order which exchanges the domain $C_3^k$ of the summary rules against the poset $\myHCC{k}{3}$. We are thus dealing with the homomorphisms contained in $\myHCCC{k}{3}{3}$. Because of
\begin{align} \label{eq_HCv1Ck1_HCvCk}
\myHCC{k-1}{3} & \; \stackrel{\eqref{isom_HPDQ_HQDP}}{\simeq} \; \Hzk,
\end{align}
the enumeration of the symmetric summary rules is equivalent to the enumeration of the homomorphisms contained in $\H( \Hzk, C_3 )$. In this section we prove
\begin{align} \label{formel_hHzk}
h( \Hzk ) & \; = \; \sum_{i=0}^k \left[ \binom{k+i}{k} - \binom{k+i}{k+1} \right] 2^{k-i}.
\end{align}
It looks like that this number is simply $\binom{ 2 k + 1}{k+1}$, but we did not make an effort to prove this equality.

\subsubsection{Recursion} \label{subsubsection_HC3Ck_rec}

\begin{figure} 
\begin{center}
\includegraphics[trim = 70 415 250 70, clip]{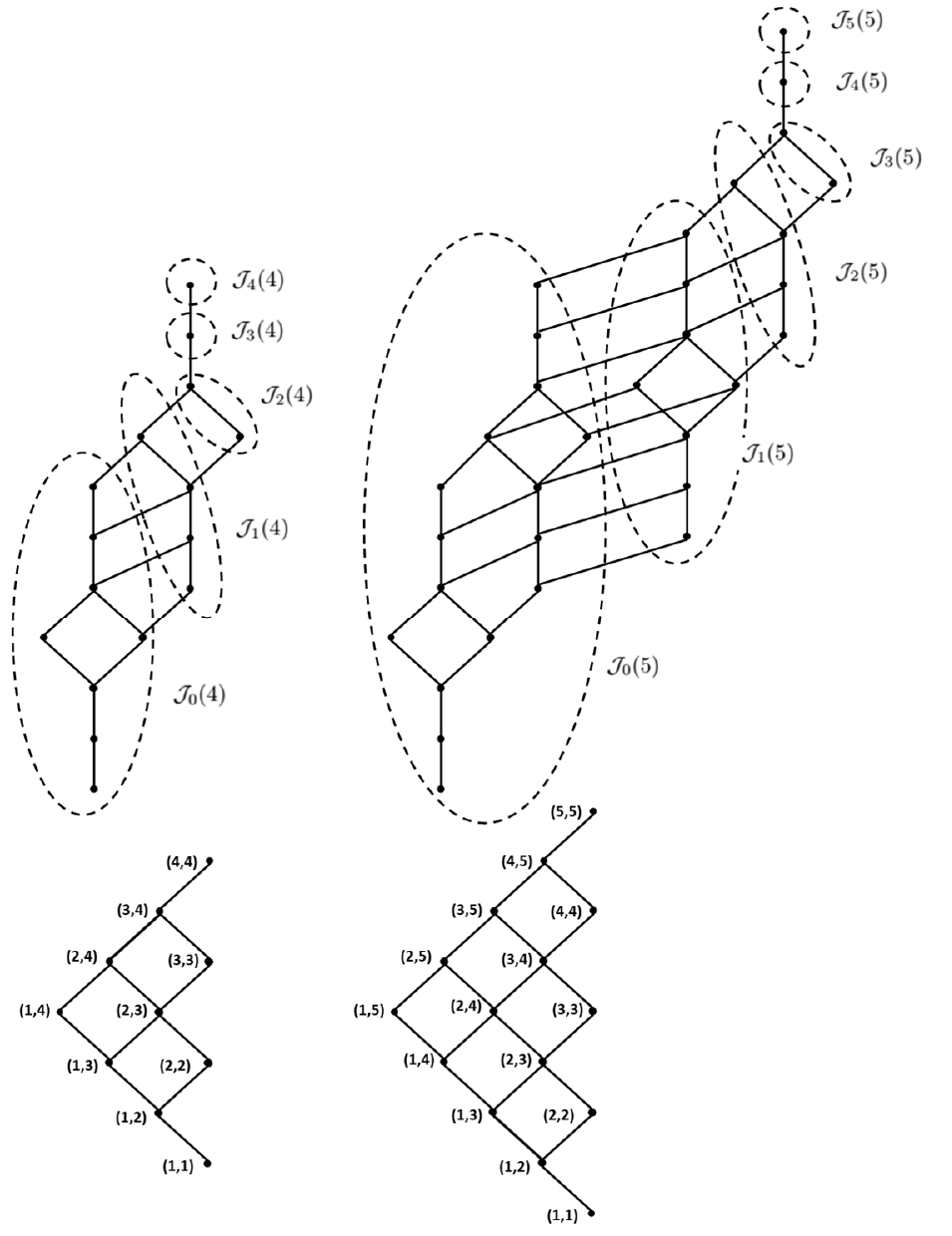}
\caption{\label{fig_HCC_DHCC} The posets $\myHCC{2}{4}$ and $\myHCC{2}{5}$ together with their down-set lattices and the respective sets $\J_j(4)$ and $\J_j(5)$.}
\end{center}
\end{figure}

We represent a homomorphism $\xi \in \Hzk$ by the pair $( \xi(1), \xi(2) )$, as shown in Figure \ref{fig_HCC_DHCC}. For $k \geq 2$, the set $\Hzk$ fits in the frame of Definition \ref{def_PQR} in the following way:
\begin{align*}
P & \equivA \myHCC{2}{k-1}, \\
Q & \equivA C_k \times \{ k \}, \\
S^- & \equivA C_{k-1} \times \{ k-1 \}, \\
S^+ & \equivA  C_{k-1} \times \{ k \}, \\
\iota(j,k) & \equivA (j,k-1) \quad \mytext{for all } j \in \Nke, \\
P^+ & \eqA \H( C_2, C_k ) \setminus \{ (k,k) \}.
\end{align*}
$Q$ is thus the top-diagonal in the Hasse diagram of $\Hzk$. Differently from Section \ref{subsec_LambdaDiamond}, we have $S^+ \not= Q$ and $P^+ \not= R$ now.

Let $k \geq 2$ be fixed. In order to unburden the notation, we use short-terms also in this section:
\begin{align*}
\J_j(k) & \equivA \J_{C_j \times \{ k \}}( \Hzk ), \\
\mytext{and} \quad a_j(k) & \equivA a_{C_j \times \{ k \}}( \Hzk )
\end{align*}
for every $j\in \Nk_0$, hence
\begin{align*}
\J_0(k) & = \mysetdescr{ D \in \Dzk }{ (1,k) \notin D }, \\
\forall j \in \Nke \mytext{:} \J_j(k) & = \mysetdescr{ D \in \Dzk }{ (j,k) \in D, (j+1,k) \notin D }, \\
\J_k(k) & = \{ \Hzk \}.
\end{align*}

Figure \ref{fig_HCC_DHCC} shows $\Hzk$ and $\Dzk$ for $k=4$ and $k=5$ together with the respective sets $\J_j(4)$ and $\J_j(5)$. The isomorphisms $\tau_T$ from Lemma \ref{lemma_beta_tau} are clearly visible. For the illustration of the isomorphism $\beta_N$ in the lemma, observe that with $N \equiv \{ (1,k) \}$, the set $\mydownarrow_{\Hzk} N$ is the bottom diagonal $(1,1), \ldots , (1,k)$ in the diagram of $\H(C_2, C_k)$, hence $\H(C_2,C_k) \smda{\Hzk} N \simeq \H(C_2,C_{k-1})$, and formula \eqref{beta_is_isom} becomes
\begin{align} \label{cup_Jjk_DHCC}
\bigcup_{j=1}^k \J_j(k) & \simeqA \D(\H(C_2,C_{k-1})),
\end{align}
as confirmed by Figure \ref{fig_HCC_DHCC} for $k = 5$.

\begin{theorem} \label{theo_card_hk}
We have
\begin{align} \label{isom_Jknull}
\J_0(k) \simeqA \D(\H(C_2,C_{k-1})) \simeq \Dzk \setminus \J_0(k),
\end{align}
hence
\begin{align} \label{cardDHC2Ck}
\# \D(\H(C_2,C_k)) & \; = \; 2^k.
\end{align}
Furthermore, $a_1(0) = 1, a_1(1) = 2$, and for every $k \geq 2$,
\begin{align} \nonumber
a_0(k) & \; = \; a_1(k) = h(\H(C_2,C_{k-1})), \\ \nonumber
\forall j \in \Nke \mytext{: } a_j(k) & \; = \; \sum_{i=j-1}^{k-1} a_{i}(k-1), \\ \nonumber
a_k(k) & \; = \; 2^k.
\end{align}
\end{theorem}
\BP The first isomorphism in \eqref{isom_Jknull} is \eqref{Jempty_isom}, and the second one is \eqref{cup_Jjk_DHCC}. With this result, \eqref{cardDHC2Ck} follows due to  $\# \D(\H(C_2,C_1)) = 2$, and \eqref{cardDHC2Ck} yields $a_k(k) = 2^k$ because of $\J_k(k) = \{ \H(C_2,C_k) \}$.

$a_0(1) = 1, a_1(1) = 2$ is easily seen. For $k \geq 2$, $a_0(k) = h( \H(C_2,C_{k-1}) )$ is due to \eqref{a_empty_R}. Let $j \in \Nke$. The right sum in \eqref{formel_a} reduces to the term $a_{ T \setminus \{ \bot \}}( P^+ \smdaPp \bot )$ with
\begin{align*}
\bot & \eqA (1,k), \\
T  & \eqA C_j \times \{ k \}, \\
T  \setminus \{ \bot \} & \simeqA C_{j-1} \times \{ k-1 \}, \\
P^+ \smdaPp \bot & \eqA 
P^+ \setminus \left( \{ 1 \} \times C_k \right) \\
& \simeqA \H(C_2, C_{k-1}) \setminus \{ (k-1,k-1) \},
\end{align*}
and because of $j - 1 < k - 1$, Lemma \ref{lemma_aTRStr_aTR} yields
\begin{align*}
a_{ T \setminus \{ \bot \}}( P^+ \smdaPp \bot ) & \eqA a_{j-1}(k-1).
\end{align*}
All together, \eqref{formel_a} delivers  $a_k(j) = \sum_{i=j-1}^{k-1} a_{k-1}(i) $; in particular, $a_1(k) = h( \H(C_2,C_{k-1}) )$.

\EP

The numbers $h( \Hzk )$ and $a_j(k)$ for $k \in \myNk{10}$ are shown in Figure \ref{fig_Tabelle_akj} in the Appendix. We have $h( \H(C_2,C_k)) = \binom{2k+1}{k}$ for all $k \in \myNk{10}$.

\subsubsection{A polynomial approach}  \label{subsubsection_HC3Ck_poly}

With our choice of $S^+$, the coefficient $a_k(k)$ was taken out of the recursion and turned out to be $2^k$. In the following definition, we use it to introduce polynomials:
\begin{definition}
For $k \in \myN$, $j \in \Nk_0$, we define polynomials $q_j^{(k)}(x)$ by setting $q_0^{(1)}(x) \equiv 1$, $q_1^{(1)}(x) \equiv x$, and, for every $k \geq 2$,
\begin{align} \nonumber 
q_0^{(k)}(x) & \equivA q_1^{(k)}(x), \\ \nonumber
j \in \Nke \mytext{: } q_j^{(k)}(x) & \equivA \sum_{i=j-1}^{k-1} q_i^{(k-1)}(x), \\ \nonumber
q_k^{(k)}(x) & \equivA x^k.
\end{align}
\end{definition}
For $j \in \Nke_0$, the degree of $q_j^{(k)}(x)$ is $k-1$, and the comparison with Theorem \ref{theo_card_hk} shows $q_0^{(k+1)}(2) = h( \Hzk )$ for all $k \in \myN$.

\begin{figure}
\begin{center}
\includegraphics[trim = 70 590 350 105, clip]{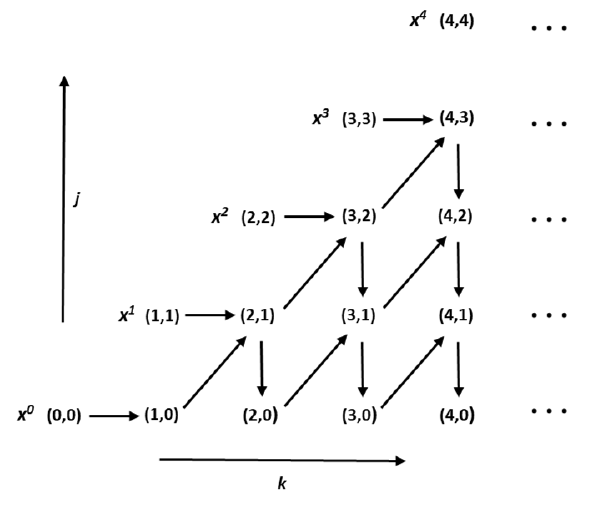}
\caption{\label{fig_HCC_Grid} The directed graph $G$ used in Lemma \ref{lemma_gridpaths}.}
\end{center}
\end{figure}

\begin{lemma} \label{lemma_gridpaths}
Let $G$ be the directed graph shown in Figure \ref{fig_HCC_Grid} with vertex set
\begin{align*}
V \equivA &  \quad \mysetdescr{ (k,j) }{ k \in \myN_0, j \in \Nk_0 }.
\end{align*}
For every $(k,j) \in V$, $i \in \Nk_0$, let $\pi_j^{(k)}(i)$ denote the number of paths in $G$ starting in $(i,i)$ and ending in $(k,j)$ (with $\pi_k^{(k)}(k) \equiv 1$). Then, for all $k \in \myN_0$,
\begin{align} \label{poly_q}
q_j^{(k)}(x) & \eqA \sum_{i=0}^{k-1} \pi_j^{(k)}(i) \cdot x^i \quad \mytext{for all } j \in \Nke_0, \\ \nonumber
\mytext{and} \quad q_k^{(k)}(x) & \eqA \sum_{i=0}^{k} \pi_k^{(k)}(i) \cdot x^i.
\end{align}
\end{lemma}

\BP The equation for $q_k^{(k)}(x)$ follows because of $\pi_k^{(k)}(i) = \delta_{ki}$ (Kronecker Delta).

Equality \eqref{poly_q} holds for $k=0$ due to $q_0^{(1)}(x) = 1 = \pi_0^{(1)}(0) \cdot x^0$. Assume that \eqref{poly_q} has been proven for $k-1 \in \myN_0$, and let $j \in \Nk$. According to the definition of $q_j^{(k)}(x)$ and the induction hypothesis, the coefficient of $x^i$ in $q_j^{(k)}(x)$ is
\begin{align} \label{sum_pi}
\sum_{\ell=j-1}^{k-1} \pi_\ell^{(k-1)}(i).
\end{align}
Every path from $(i,i)$ to $(k,j)$ has to run over one of the vertices $(k-1,\ell)$ with $j-1 \leq \ell \leq k$. For $j-1 \leq \ell \leq k$, let $\fP(\ell)$ be the set of paths in $G$ from $(i,i)$ to $(k-1,\ell)$. Each path in $\fP(\ell)$ can be extended to a path from $(i,i)$ to $(k,j)$ by adding a step diagonally right upwards to $(k,\ell+1)$ followed by $ \ell + 1 - j$ steps downwards. Defining for all $j-1 \leq \ell \leq k$,
\begin{align*}
\fP'(\ell) & \equivA \mysetdescr{ \fp (k,\ell+1) (k,\ell) \cdots (k,j) }{ \fp \in \fP(\ell) },
\end{align*}
the sets $\fP'(\ell)$ are pairwise disjoint with $\# \fP'(\ell) = \# \fP(\ell) = \pi_\ell^{(k-1)}(i)$. The set $\cup_{\ell=j}^{k-1} \fP'(\ell)$ contains the paths from $(i,i)$ to $(k,j)$ running over $(k,j+1)$, and the set $\fP'(j-1)$ contains the paths from $(i,i)$ to $(k,j)$ running over $(k-1,j-1)$, and \eqref{sum_pi} is the number of paths from $(i,i)$ to $(k,j)$.

The polynomial $q_0^{(k)}(x)$ remains. Because of $k > 0$, every path from $(i,i)$ to $(k,0)$ must run over $(k,1)$. The number $\pi_0^{(k)}(i) = \pi_1^{(k)}(i)$ is thus the number of paths in $G$ from $(i,i)$ to $(k,0)$.

\EP

Equation \eqref{formel_hHzk} is a direct consequence of the following theorem:

\begin{theorem} \label{theo_coeffs_bkj}
For all $k \in \myN_0$,
\begin{align} \label{formula_qkez}
q_0^{(k+1)}(x) & \; = \; \sum_{i=0}^k \left[ \binom{k+i}{k} - \binom{k+i}{k+1} \right] x^{k-i}.
\end{align}
\end{theorem}
\BP Let $k \in \myN_0$. According to Lemma \ref{lemma_gridpaths}, we have
\begin{align*}
q_0^{(k+1)}(x) & \; = \; \sum_{i=0}^k \pi_0^{(k+1)}(i) \cdot x^i,
\end{align*}
where $\pi_0^{(k+1)}(i)$ is the number of paths in $G$ from $(i,i)$ to $(k+1,0)$. It is also the number of paths in $G$ from $(i+1,i)$ to $(k+1,0)$. Each of these paths we can step backwards from $(k+1,0)$ to $(i+1,i)$. Such a reversed path can be described by an unique sequence of $k$ letters $U$ (for steps upwards in $G$) and $k-i$ letters $D$ (for diagonal steps downwards) in which the $\ell$\textsuperscript{th} occurrence of $D$ is preceded by at least $\ell$ occurences of $U$. Let $u_1, \ldots , u_k$ be the indexes of the $U$-letters in such a sequence, and $d_1, \ldots , d_{k-i}$ the indexes of the $D$-letters. Writing $u_1, \ldots , u_k$ into the upper row of a Ferrers diagram of type $[k,k-i]$ and $d_1, \ldots , d_{k-i}$ into the lower one yields a standard Young-tableau of type $[k,k-i]$. This mapping between reversed $(i+1,i)$-$(k+1,0)$-paths in $G$ and standard Young-tableaux of type $[k,k-i]$ is bijective, and application of the hook-length-formula of Frame, Robinson, and Thrall yields
\begin{align*}
\pi_0^{(k+1)}(i) & \eqA \frac{( 2k-i )!}{ \frac{(k+1)!}{i+1} \cdot (k-i)!} \eqA \binom{2k-i}{k} - \binom{2k-i}{k+1}.
\end{align*}

\EP

The polynomial coefficients $\binom{k+i}{k} - \binom{k+i}{k+1}$ are shown in Figure \ref{fig_Tabelle_PolCoeffs} in the Appendix for $k \in \myNk{10}$ and $i \in \Nk_0$.
\newline

\noindent {\bf Acknowledgement:} I am grateful to Lawrence C.\ Rafsky for making me enthusiastic about the combinatorics of $\H(R,C_3)$.

\newpage

\section{Appendix}

\begin{figure}[ht]
\begin{center}
\includegraphics[trim = 70 675 210 70, clip]{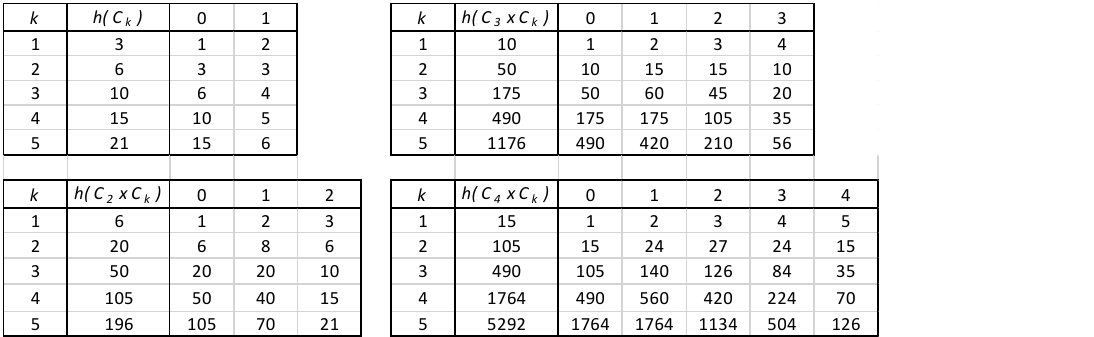}
\caption{\label{fig_Tabelle_CnCk} The coefficients $h( \CCnk )$ and $a_j(\CCnk)$ for $n \in \myNk{4}$ and $k \in \myNk{5}$.}
\end{center}
\end{figure}

\begin{figure}[ht]
\begin{center}
\includegraphics[trim = 70 580 210 70, clip]{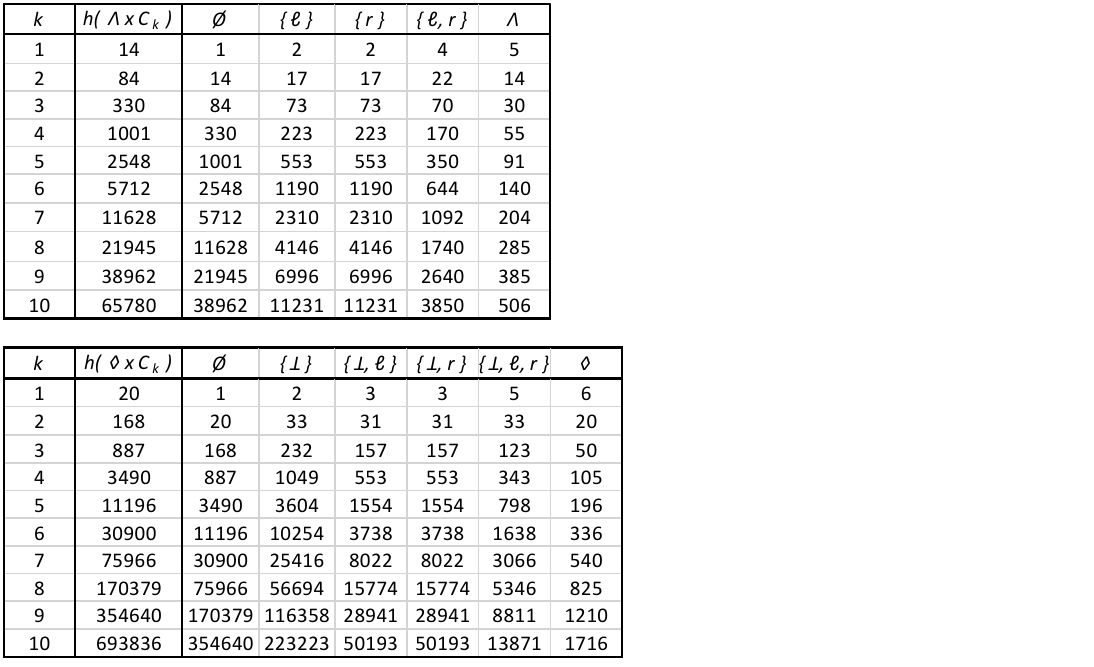}
\caption{\label{fig_LambdaLozenge} $h( \LCk )$, $a_T( \LCk )$, $h( \LOCk )$, and $a_T( \LOCk )$ for $k \in \myNk{10}$. We have $a_\Lambda(\LCk) = \# \D( \LCk )$ and $a_\lozenge(\LOCk) = \# \D( \LOCk )$.}
\end{center}
\end{figure}

\clearpage
\begin{figure}[ht]
\begin{center}
\includegraphics[trim = 70 680 210 70, clip]{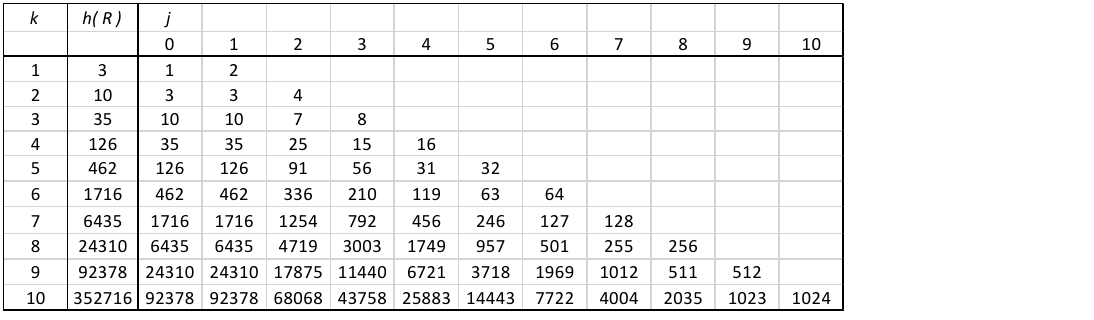}
\caption{\label{fig_Tabelle_akj} The numbers $h( \Hzk )$ and $a_j(k)$ for $k \in \myNk{10}$.}
\end{center}
\end{figure}

\begin{figure}[ht]
\begin{center}
\includegraphics[trim = 70 680 210 70, clip]{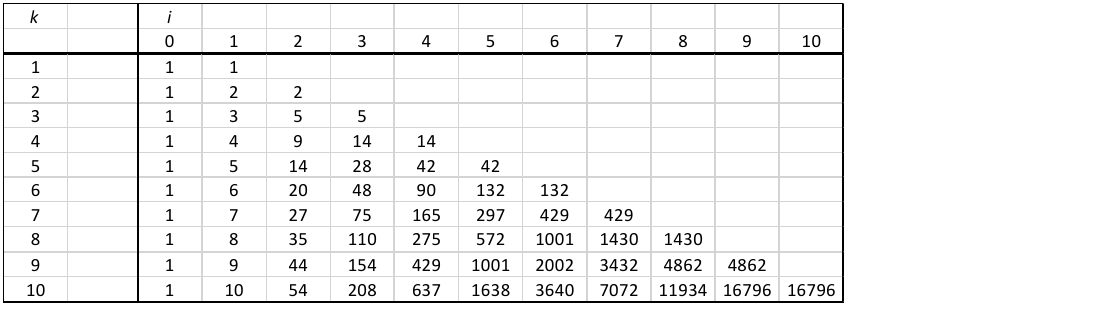}
\caption{\label{fig_Tabelle_PolCoeffs} The polynomial coefficients $\binom{k+i}{k} - \binom{k+i}{k+1}$ for $k \in \myNk{10}$.}
\end{center}
\end{figure}


\begin{thebibliography}{xx}
{\small
\bibitem{Aigner_1975} M. Aigner: {\em Kombinatorik I. Grundlagen und Z\"{a}hltheorie}. Springer-Verlag Berlin, Heidelberg, New York (1975).


\bibitem{Brightwell_Winkler_1991} G. Brightwell and P. Winkler: Counting linear extensions. {\em Order} {\bf 8 } (1991), 225--242.

\bibitem{Browning_etal_2018} T. Browning, M. Hopkins, and Z. Kelley: Doppelgangers: the ur-operation and posets of bounded height. arXiv:1710.10407v2 [math.CO] (2018).

\bibitem{aCampo_2018} F. a Campo: Relations between powers of Dedekind numbers and exponential sums related to them. {\em J. Int. Seq.}  {\bf 21} (2018), Article 18.4.4.

\bibitem{aCampo_2019} F. a Campo: A framework for the systematic determination of the posets on $n$ points with at least $\tau \cdot 2^n$ downsets. {\em Order} {\bf 36} (2019), 119--157. Published Online May 29, 2018, https://doi.org/10.1007/s11083-018-9459-2.

\bibitem{aCampo_Erne_2019} F. a Campo and M. Ern\'{e}: Exponential functions of finite posets and the number of extensions with a fixed set of minimal points. {\em J.\ Comb.\ Math.\ and Comb.\ Calc.\ } {\bf 110} (2019), 125--156.


\bibitem{Edelman_Klingsberg_1982} P. H. Edelman and  P. Klingsberg: The subposet lattice and the order polynomial. {\em Europ. J. Combinatorics} {\bf 3} (1982), 341--346.

\bibitem{Hamaker_etal_2020} Z. Hamaker, R. Patrias, O. Pechenik, and N. Williams: Doppelg\"{a}ngers: Bijections of plane partitions. {\em International Mathematics Research Notices} {\bf 2020} (2020), 487--540. Pre-published in 2016.

\bibitem{Hopkins_2020} S. Hopkins: Order polynomial product formulas and poset dynamics. arXiv:2006.01568v3 [math.CO] (2020).

\bibitem{Jochemko_2014} K. Jochemko: Order polynomials and P\'{o}lya's enumeration theorem. {\em The electronic journal of combinatorics} {\bf 21} (2014), $\#$P2.52.

\bibitem{Kovalerchuk_1995} B. Kovalerchuk, E. Triantaphyllou, and E. Vityaev: Monotone boolean function learning techniques integrated with user interaction. In {\em Proceedings of Workshop ``Learning from Examples vs.\ Programming by Demonstration'', 12\textsuperscript{th} International Conference on Machine Learning}, Tahoe City, CA, (1995), 41--48.

\bibitem{Liggins_Nebrich_2000} M. E. Liggins II and M. A. Nebrich: Adaptive multi-image decision-fusion. In I. Kadar (ed.): {\em Signal Processing, Sensor Fusion, and Target Recognition IX}, Proceedings of SPIE {\bf 4052} (2000), 218--228.

\bibitem{Polikar_2006} R. Polikar: Ensemble based systems in decision making. {\em IEEE Circuit Syst. Mag.} {\bf 6-3} (2006), 21--45.

\bibitem{Stanley_1970} R. P. Stanley A chromatic-like polynomial for ordered sets. {\em Proceedings of the
Second Chapel Hill Conference on Combinatorial Mathematics and its Applications} (1970), 421--427.

\bibitem{Stanley_1971_Diss} R. P. Stanley: Ordered Structures and Partitions. PhD thesis, Harvard (1971).

\bibitem{Stanley_1972} R. P. Stanley: Ordered Structures and Partitions. Memoirs of the AMS (1972).

\bibitem{Stanley_1984} R. P. Stanley: Two poset polytopes. {\em Discrete Comput.\ Geom.} {\bf 1} (1986), 9--23.

\bibitem{Stanley_2012} R. P. Stanley: {\em Enumerative Combinatorics, Volume I}. Cambridge University Press, 2\textsuperscript{nd} edition (2012).


\bibitem{Thomas_2003} H. Thomas: Order-preserving maps from a poset to a chain, the order polytope, and the Todd class of the associated toric variety. {\em Europ. J. Combinatorics} {\bf 24} (2003), 809--814.

\bibitem{Wagner_1992} D. G. Wagner: Enumeration of functions from posets to chains. {\em Europ. J. Combinatorics} {\bf 13} (1992), 313--324.
}
\end{thebibliography}
\end{document}